\documentclass[a4paper, 11pt]{article}

\usepackage[latin1]{inputenc}
\usepackage[english]{babel}
\usepackage{graphicx}
\usepackage{amssymb}
\usepackage{amsmath}
\usepackage{amsfonts}
\usepackage{amsthm}
\usepackage{url}
\usepackage[authoryear]{natbib}

\newtheorem{theorem}{Theorem}

\newtheorem{algorithm}{Algorithm}

\numberwithin{equation}{section}
\DeclareMathOperator{\argmin}{argmin}
\usepackage{algorithm}
\usepackage{algorithmic}
\def \matlab    {MATLAB$^{\textrm{\tiny \textregistered}}$}
\DeclareMathOperator{\update}{update}

\title{Accelerated Multiplicative Updates and Hierarchical ALS Algorithms for Nonnegative Matrix Factorization}

\author{\normalsize Nicolas Gillis${}^1$ and Fran\c{c}ois Glineur${}^2$}

\date{}

\begin{document}

\maketitle

\begin{abstract}
Nonnegative matrix factorization (NMF) is a data analysis technique used in a great variety of applications such as text mining, image processing, hyperspectral data analysis, computational biology, and clustering. In this paper, we consider two well-known algorithms designed to solve NMF problems, namely the multiplicative updates of Lee and Seung and the hierarchical alternating least squares of Cichocki et al. We propose a simple way to significantly accelerate these schemes, based on a careful analysis of the computational cost needed at each iteration, while preserving their convergence properties. This acceleration technique can also be applied to other algorithms, which we illustrate on the projected gradient method of Lin. The efficiency of the accelerated algorithms is empirically demonstrated on image and text datasets, and compares favorably with a state-of-the-art alternating nonnegative least squares algorithm. \bigskip

\noindent {\bf Keywords:} nonnegative matrix factorization,  algorithms, multiplicative updates, hierarchical alternating least squares. 
\end{abstract}

\footnotetext[1] {University of Waterloo, Department of Combinatorics and Optimization, Waterloo, Ontario N2L 3G1, Canada. E-mail: ngillis@uwaterloo.ca. This work was carried out when the author was a Research fellow of the Fonds de la Recherche Scientifique (F.R.S.-FNRS) at Universit\'e catholique de Louvain.} 

\footnotetext[2] {Universit\'e catholique de Louvain, CORE and ICTEAM Institute, B-1348 Louvain-la-Neuve, Belgium. 
E-mail: francois.\mbox{glineur}@uclouvain.be.  This text presents research results of the Belgian Program on Interuniversity Poles of Attraction initiated by the Belgian State, Prime Minister's Office, Science Policy Programming. The scientific responsibility is assumed by the authors.}

\section{Introduction} \label{intro}

Nonnegative matrix factorization (NMF) consists in approximating a nonnegative matrix $M$ as a low-rank product of two nonnegative matrices $W$ and $H$, i.e., given a matrix $M \in \mathbb{R}^{m \times n}_+$ and an integer $r < \min\{m,n\}$, find two matrices $W\in \mathbb{R}^{m \times r}_+$ and $H\in\mathbb{R}^{r \times n}_+$ such that $WH \approx M$.

With a nonnegative input data matrix $M$, nonnegativity constraints on the factors $W$ and $H$ are well-known to lead to low-rank decompositions with better interpretation in many applications such as  text mining \citep{SBPP06}, image processing \citep{LS1}, hyperspectral data analysis \citep{PPP06}, computational biology \citep{Dev}, and clustering \citep{DHS05}. Unfortunately, imposing these constraints is also known to render the problem computationally difficult \citep{V09}. 

Since an exact low-rank representation of the input matrix does not exist in general, the quality of the approximation is measured by some criterion, typically the sum of the squares of the errors on the entries, which leads to the following minimization problem:
\begin{equation}
\min_{W \in \mathbb{R}^{m \times r}, H \in \mathbb{R}^{r \times n}} 
||M-WH||_F^2 \tag{NMF} \quad \text{ such that } \quad
 W \geq 0 \text{ and } H \geq 0,  \label{NMF} 
 \end{equation} 
 where $||A||_F = (\sum_{i,j} A_{ij}^2)^\frac12$ denotes the Frobenius norm of matrix $A$. 
 Most NMF algorithms are iterative, and exploit the fact that \eqref{NMF} reduces to an efficiently solvable convex nonnegative least squares problem (NNLS) when one of the factors $W$ or $H$ is fixed. Actually, it seems that nearly all algorithms proposed for NMF adhere to the following general framework 
 \begin{enumerate} 
 \item[(0)] Select initial matrices $(W^{(0)}, H^{(0)})$ (e.g., randomly). Then for $k = 0, 1, 2, \dots$, do
 \item[(a)] Fix $H^{(k)}$ and  find $W^{(k+1)} \geq 0$  such that $||M-W^{(k+1)}H^{(k)}||_F^2 < ||M-W^{(k)}H^{(k)}||_F^2$. 
 \item[(b)] Fix $W^{(k+1)}$ and  find $H^{(k+1)} \geq 0$ such that $||M-W^{(k+1)}H^{(k+1)}||_F^2 < ||M-W^{(k+1)}H^{(k)}||_F^2$.
 \end{enumerate} 
More precisely, at each iteration, one of the two factors is fixed and the other is updated in such a way that the objective function is reduced, which amounts to a two-block coordinate descent method. Notice that the role of matrices $W$ and $H$ is perfectly symmetric: if one transposes input matrix $M$, the new matrix $M^T$ has to be approximated by a product $H^T W^T$, so that any formula designed to update the first factor in this product directly translates into an update for the second factor in the original problem. Formally, if the update performed in  step (a) is described by $W^{(k+1)} = \update(M,W^{(k)},H^{(k)})$, an algorithm preserving symmetry will update the factor in step (b) according to $H^{(k+1)} = \update(M^T,H^{(k)}{}^T,W^{(k+1)}{}^T)^T$. In this paper, we only consider such symmetrical algorithms, and focus on the update of matrix $W$. \

This update can be carried out in many different ways: the most natural possibility is to compute an optimal solution for the NNLS subproblem, which leads to a class of algorithms called alternating nonnegative least squares (ANLS), see, e.g., \cite{KP08}. 
However, this computation, which can be performed with active-set-like methods \citep{KP08, KP09}, is relatively costly. Therefore, since an optimal solution for the NNLS problem corresponding to one factor is not required before the update of the other factor is performed, several algorithms only compute an approximate solution of the NNLS subproblem, sometimes very roughly, but with a cheaper computational cost, leading to an inexact two-block coordinate descent scheme. We now present two such procedures: the multiplicative updates of Lee and Seung and the hierarchical alternating least squares of Cichocki et al.

In their seminal papers, \cite{LS1, LS2} introduce the multiplicative updates: 
\begin{equation} 
W^{(k+1)} = \text{MU}(M,W^{(k)},H^{(k)}) = W^{(k)} \circ \frac{[M {H^{(k)}}^T]}{[W^{(k)} H^{(k)} {H^{(k)}}^T]}, 
\nonumber 
\end{equation}
where $\circ$ (resp.\@ {\footnotesize{$\frac{[\,.\,]}{[\,.\,]}$}}) denotes the component-wise product (resp.\@ division) of matrices, and prove that each update monotonically decreases the Frobenius norm of the error $||M-WH||_F$, i.e., satisfies the description of steps (a) and (b). This technique was actually originally proposed by \cite{DM86} to solve NNLS problems. 
The popularity of this algorithm came along with the popularity of NMF and many authors have studied or used this algorithm or variants to compute NMF's, see, e.g., \cite{Ber, CAZP09} and the references therein. In particular, the \matlab \, Statistics Toolbox implements this method. 

However, MU have been observed to converge relatively slowly,  especially when dealing with dense matrices $M$, see \cite{HHNP09, GG08} and the references therein, and many other algorithms have been subsequently  introduced which perform better in most situations. For example, \cite{Cic, Cic4} and, independently, several other authors \citep{diep, GG08, LZ09} proposed a technique called hierarchical alternating least squares (HALS)\footnote{\cite{diep} refers to HALS as rank-one residue iteration (RRI), and \cite{LZ09} as FastNMF.}, which successively updates each column of $W$ 
with an optimal and easy to compute closed-form solution. In fact, 
when fixing all variables but a single column $W_{:p}$ of $W$, 
the problem reduces to 
\[
\min_{W_{:p} \geq 0} ||M-WH||_F^2 
= ||(M-\sum_{l\neq p}W_{:l}H_{l:}) - W_{:p}H_{p:}||_F^2 
= \sum_{i = 1}^m ||(M_{i:}-\sum_{l\neq p}W_{il}H_{l:}) - W_{ip}H_{p:}||_F^2. 
\]
Because each row of $W$ only affects the corresponding row of the product $W H$, this problem can be further decoupled into $m$ independent quadratic programs in one variable $W_{ip}$, corresponding to the $i^{\text{th}}$ row of $M$. The optimal solution $W_{ip}^*$ of these subproblems can be easily written in closed-form 
\begin{eqnarray*}
W_{ip}^* 
& = & \max\Big(0,\frac{(M_{i:}-\sum_{l\neq p}W_{il}H_{l:})H_{p:}^T}{H_{p:}H_{p:}^T} \Big) \\ 
& = & \max\Big(0, \frac{M_{i:}H_{p:}^T-\sum_{l \neq p}W_{il}H_{l:}H_{p:}^T}{H_{p:}H_{p:}^T} \Big), \; 1 \leq i \leq m.
\end{eqnarray*}
Hence HALS updates successively the columns of $W$, so that 
$W^{(k+1)} = \text{HALS}(M,W^{(k)},H^{(k)})$ can be computed in the following way:   
%
\[
W_{:p}^{(k+1)} = \max\Big({0}, \frac{A_{:p} - \sum_{l = 1}^{p-1} W_{:l}^{(k+1)} B_{lp}  - \sum_{l = p+1}^{r} W_{:l}^{(k)} B_{lp}}{B_{pp}} \Big), 
\] 
$\text{successively for } p = 1, 2, \dots, r$, where $A = M{H^{(k)}}^T$ and $B = H^{(k)}{H^{(k)}}^T$. 
This amounts to approximately solving each NNLS subproblem in $W$ with a single complete round of an exact block-coordinate descent method with $r$ blocks of $m$ variables corresponding to the columns of $W$ (notice that any other ordering for the update of the columns of $W$ is also possible).  


Other approaches based on iterative methods to solve the NNLS subproblems include projected gradient descent \citep{Lin} or Newton-like methods~\citep{DKS07, CZA06}; see also \cite{CAZP09} and the  references therein. \\

We first analyze in Section~\ref{ccost} the computational cost needed to update the factors $W$ in MU and HALS, then make several simple observations leading  in Section \ref{accelalgo} to the design of accelerated versions of these algorithms. These improvements can in principle be applied to any two-block coordinate descent NMF algorithm, as demonstrated in Section~\ref{appllin} on the projected gradient method of \cite{Lin}. We mainly focus on MU, because it is by far the most popular NMF algorithm, and on HALS, because it is very efficient in practice. 
Section~\ref{conv} studies convergence of the accelerated variants to stationary points, and shows that they preserve the properties of the original schemes.
In Section~\ref{ne}, we experimentally demonstrate a significant acceleration in convergence on several image and text datasets, with a comparison with the state-of-the-art ANLS algorithm of \cite{KP09}. 


\section{Analysis of the Computational Cost of Factor Updates} \label{ccost}



In order to make our analysis valid for both dense and sparse input matrices,  
let us introduce a parameter $K$ denoting the number of nonzero entries in matrix $M$ ($K = mn$ when $M$ is dense).  Factors $W$ and $H$ are typically stored as dense matrices throughout the execution of the algorithms. We assume that NMF achieves compression, which is often a requirement in practice. This means that storing $W$ and $H$ must be cheaper than storing $M$:  roughly speaking, the number of entries in $W$ and $H$ must be smaller than the number of nonzero entries in $M$, i.e., $r(m+n) \le K$. 

Descriptions of Algorithms~\ref{MUalgo} and \ref{HALSalgo} 
below provide separate estimates for the number of floating point operations (flops) in each matrix product computation needed to update factor $W$ in MU and HALS
.  
One can check that the proposed organization of the different matrix computations (and, in particular, the ordering of the matrix products) minimizes the total computational cost (for example, starting the computation of the MU denominator $WHH^T$ with the product $WH$ is clearly worse than with $HH^T$). 
\algsetup{indent=2em} 
\begin{algorithm}
\caption{MU update for $W^{(k)}$} \label{MUalgo}
\begin{algorithmic}[1]
\STATE $A = M {H^{(k)}}^T$;  \hspace{1.13cm}  
$\quad\rightarrow\quad 2Kr$ flops
\STATE $B = H^{(k)} {H^{(k)}}^T$; \hspace{0.78cm}   
$\quad\rightarrow\quad 2nr^2$ flops 
\STATE $C = W^{(k)} B$;    \hspace{1.42cm}   
$\quad\rightarrow\quad 2mr^2$ flops 
\STATE $W^{(k+1)} =  W^{(k)} \circ \frac{[A]}{[C]}$;     
$\quad \, \rightarrow\quad 2mr$ flops \vspace{0.1cm} \\
 \% Total:  $r (2K+2nr+2mr+2m)$ flops 
\end{algorithmic}
\end{algorithm}
\algsetup{indent=2em}
\begin{algorithm}[ht!]
\caption{HALS update for $W^{(k)}$}  \label{HALSalgo}
\begin{algorithmic}[1]
\STATE $A = M {H^{(k)}}^T$;  \hspace{5.5cm}  
$\quad\rightarrow\quad 2Kr$ flops
\STATE $B = H^{(k)} {H^{(k)}}^T$; \hspace{5.15cm}   
$\quad\rightarrow\quad 2nr^2$ flops 
\FOR {$i = 1, 2, \dots, r$} 
\STATE $C_{:k} = \sum_{l = 1}^{p-1} W_{:l}^{(k+1)} B_{lk}  + \sum_{l = p+1}^{r} W_{:l}^{(k)} B_{lk}$;  			
 $\; \;  \rightarrow\quad 2m(r-1)$ flops 
\STATE $W_{:k} = \max\Big(0,\frac{A_{:k}-C_{:k}}{B_{kk}}\Big)$; \hspace{3cm}  
$\; \; \, \rightarrow\quad 3m$ flops 
\ENDFOR \vspace{0.1cm} \\
 \% Total:  $r (2K+2nr+2mr+m)$ flops 
\end{algorithmic}
\end{algorithm}

MU and HALS possess almost exactly the same computational cost (the difference being a typically negligible $mr$ flops). It is particularly interesting to observe that 
\begin{enumerate}
\item Steps 1.\@ and 2.\@ in both algorithms are identical and do not depend on the matrix $W^{(k)}$; 
\item Recalling our assumption $K \geq r(m+n)$, computation of $M{H^{(k)}}^T$ (step 1.\@) is the most expensive among all steps. 
\end{enumerate}
Therefore, this time-consuming step should be performed sparingly, and we should take full advantage of having computed the relatively expensive $M{H^{(k)}}^T$ and ${H^{(k)}}{H^{(k)}}^T$ matrix products. This can be done by updating $W^{(k)}$ several times before the next update of $H^{(k)}$, i.e., by repeating steps 3.\@ and 4.\@ in MU (resp.\@ steps 3.\@ to 6.\@ in HALS) several times after the computation of matrices $M{H^{(k)}}^T$ and ${H^{(k)}}{H^{(k)}}^T$. In this fashion, better solutions of the corresponding NNLS subproblems will be obtained at a relatively cheap additional cost.

The original MU and HALS algorithms do not take advantage of this fact, and alternatively update matrices $W$ and $H$ only once per (outer) iteration.  
An important question for us is now: 
how many times should we update $W$ per outer iteration?, i.e., how many inner iterations of MU and HALS should we perform? 
This is the topic of the next section.

\section{Stopping Criterion for the Inner Iterations} \label{accelalgo}

In this section, we discuss two different strategies for choosing the number of inner iterations: the first uses a fixed number of inner iterations determined by the flop counts, while the second is based on a dynamic stopping criterion that checks the difference between two consecutive iterates. The first approach is shown empirically to work better. We also describe a third hybrid strategy that provides a further small improvement in performance.

\subsection{Fixed Number of Inner Iterations}  \label{fnit}

Let us focus on the MU algorithm (a completely similar analysis holds for HALS, as both methods differ only by a negligible number of flops). Based on the flops counts, we estimate how expensive the first inner update of $W$ would be relatively to the next ones (all performed while keeping $H$ fixed), which is given by the following factor $\rho_W$  
(the corresponding value for $H$ will be denoted by $\rho_H$) 
\[
\rho_W 
= \frac{2Kr + 2nr^2 + 2mr^2 + 2 mr}{2mr^2 + 2 mr} = 1 + \frac{K + nr}{mr + m}. 
\quad \Big(\rho_H = 1 + \frac{K + mr}{nr + n}\Big). 
\]
Values of $\rho_W$ and $\rho_H$ for several datasets are given in Section~\ref{ne}, see Tables~\ref{dip} and \ref{dtm}. 

Notice that for $K \geq r(m+n)$, we have $\rho_W \geq 2\frac{r}{r+1}$ so that the first inner update of $W$ is at least about twice as expensive as the subsequent ones. 
For a dense matrix, $K$ is equal to $mn$ and we actually have that $\rho_W = 1+\frac{n(m + r)}{m (r + 1)} \geq 1+ \frac{n}{r+1}$, which is typically quite large since $n$ is often much greater than $r$. 
This means for example that, in our accelerated scheme, $W$ could be updated about $1+\rho_W$ 
times for the same computational cost as two independent updates of $W$ in the original MU.  

A simple and natural choice consists in performing inner updates of $W$ and $H$ 
a fixed number of times, depending on the values of $\rho_W$ and $\rho_H$. 
Let us introduce a parameter  $\alpha \geq 0$ such that $W$ is updated $(1 + \alpha \rho_W)$~times before the next update of $H$, and $H$ is updated $(1 + \alpha \rho_H)$~times before the next update of $W$. Let us also  denote the corresponding algorithm MU$_{\alpha}$ (MU$_{0}$ reduces to the original MU). 
Therefore, performing $(1 + \alpha \rho_W)$ inner updates of $W$ in MU$_{\alpha}$ has approximately  the same computational cost as performing $(1+\alpha)$ updates of $W$ in MU$_{0}$.  


In order to find an appropriate value for parameter $\alpha$, we have performed some preliminary tests on image and text datasets. First, let us denote $e(t)$ the Frobenius norm of the error $||M-WH||_F$ achieved by an algorithm within time $t$,  and define 
\begin{equation} \label{Evol} 
E(t) = \frac{e(t)-e_{\min}}{e(0)-e_{\min}}, 
\end{equation}
where $e(0)$ is the error of the initial iterate $(W^{(0)}, H^{(0)})$, and $e_{\min}$ is the smallest error observed among all algorithms across all initializations. 
Quantity $E(t)$ is therefore a normalized measure of the improvement of the objective function (relative to the initial gap) with respect to time; we have $0 \leq E(t) \leq 1$ for monotonically decreasing algorithms (such as MU and HALS). The advantage of $E(t)$ over $e(t)$ is that one can meaningfully take the average over several runs involving different initializations and datasets, and display the average behavior of a given algorithm.

\begin{figure*}[ht!]
\begin{center}
\includegraphics[width=\textwidth]{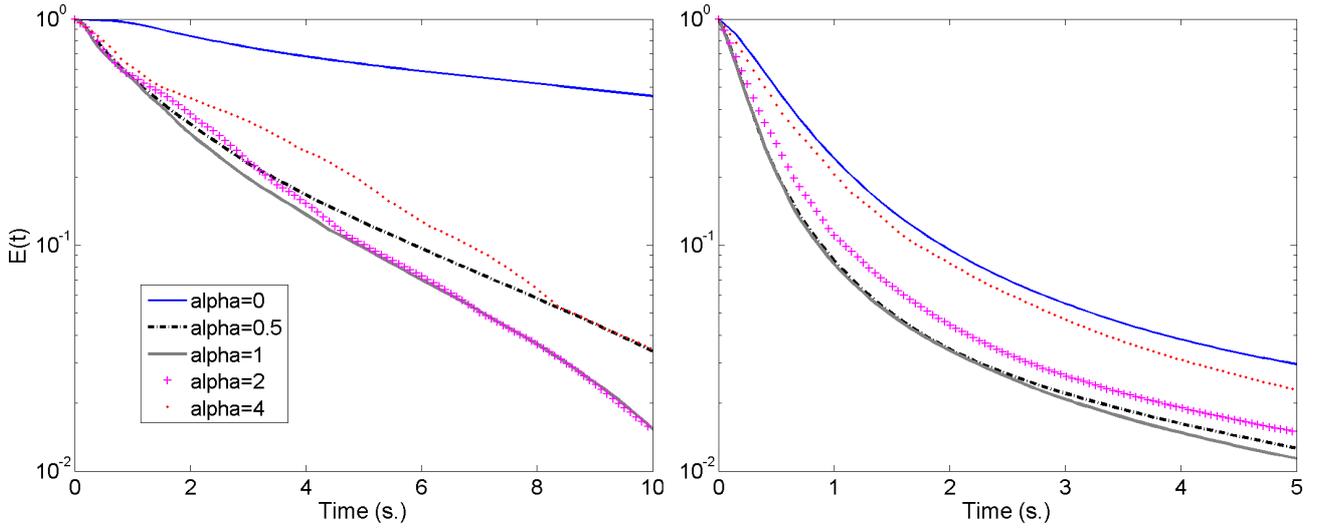}
\caption{Average of functions $E(t)$ for MU using different values of $\alpha$: (left) dense matrices, (right) sparse matrices, computed over 4 image datasets and 6 text datasets, using two different values for the rank for each dataset and 10 random initializations, see Section~\ref{ne}.} 
\label{mualpha}
\end{center}
\end{figure*}

Figure~\ref{mualpha} displays the average of this function $E(t)$ for dense (on the left) and sparse (on the right) matrices using the datasets described in Section~\ref{ne} for five values of $\alpha=0,0.5,1,2,4$.  
We observe that the original MU algorithm ($\alpha = 0$) converges significantly less rapidly than all the other tested variants (especially in the dense case). The best value for parameter $\alpha$ seems to be $1$.  

\begin{figure*}[ht!]
\begin{center}
\includegraphics[width=\textwidth]{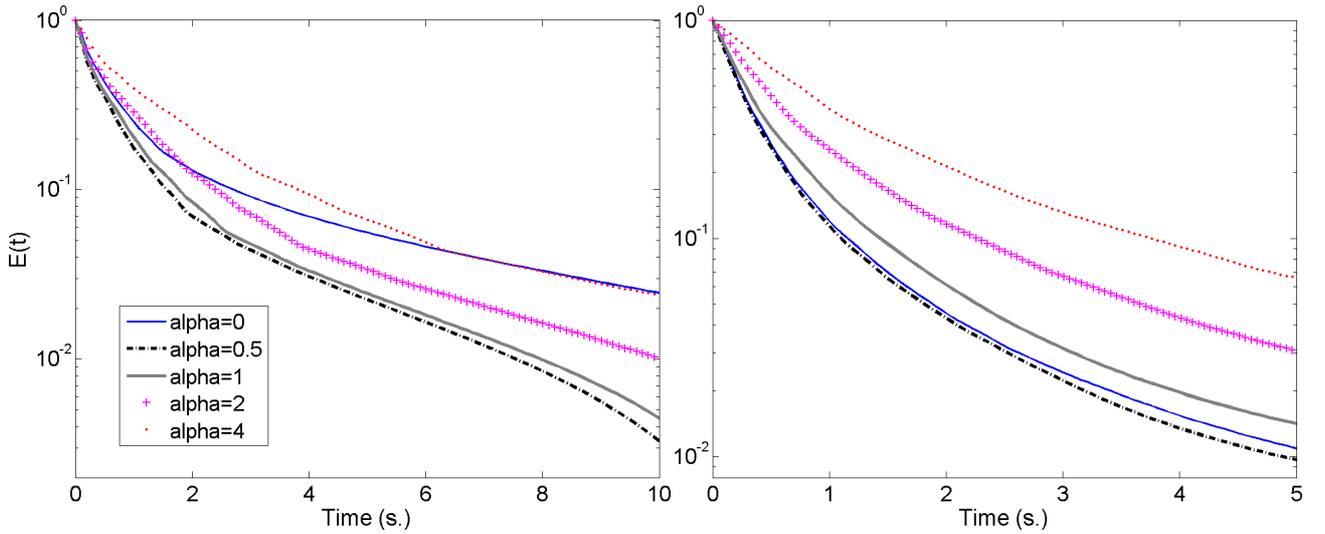}
\caption{Average of functions $E(t)$ for HALS using different values of $\alpha$: (left) dense matrices, (right) sparse matrices. Same settings as in Figure~\ref{mualpha}.} 
\label{halsalpha}
\end{center}
\end{figure*}

Figure~\ref{halsalpha} displays the same computational experiments for 
HALS\footnote{Because HALS involves a loop over the columns of $W$ and rows of $H$, we observed that an update of HALS is noticeably slower than an update of MU when using \matlab (especially for $r \gg 1$), despite the quasi-equivalent theoretical computational cost.  
Therefore, to obtain fair results, we adjusted $\rho_W$ and $\rho_H$ 
by measuring directly the ratio between time spent for the first update and the next one, using the \emph{cputime} function of \matlab.}. 
HALS with $\alpha=0.5$ performs the best. For sparse matrices, the improvement is harder to discern (but still present); an explanation for that fact will be given at the end of Section~\ref{hybrid}. 

%

\subsection{Dynamic Stopping Criterion for Inner Iterations}   \label{dcit}

In the previous section, a fixed number of inner iterations is performed. One could instead consider switching dynamically from one factor to the other based on an appropriate criterion. 
For example, it is possible to use the norm of the projected gradient as proposed by \cite{Lin}. A simpler and cheaper possibility is to rely solely on the norm of the difference between two iterates.  
Noting $W^{(k,l)}$ the iterate after $l$ updates of $W^{(k)}$ (while $H^{(k)}$ is being kept fixed), we stop inner iterations as soon as  
\begin{equation} \label{stopcrit}
||W^{(k,l+1)}-W^{(k,l)}||_F \leq \epsilon ||W^{(k,1)}-W^{(k,0)}||_F,  
\end{equation}
i.e., as soon as the improvement of the last update becomes negligible compared to the one obtained with the first update, 
but without any a priori fixed maximal number of inner iterations. 

\begin{figure*}[ht!]
\begin{center}
\includegraphics[width=\textwidth]{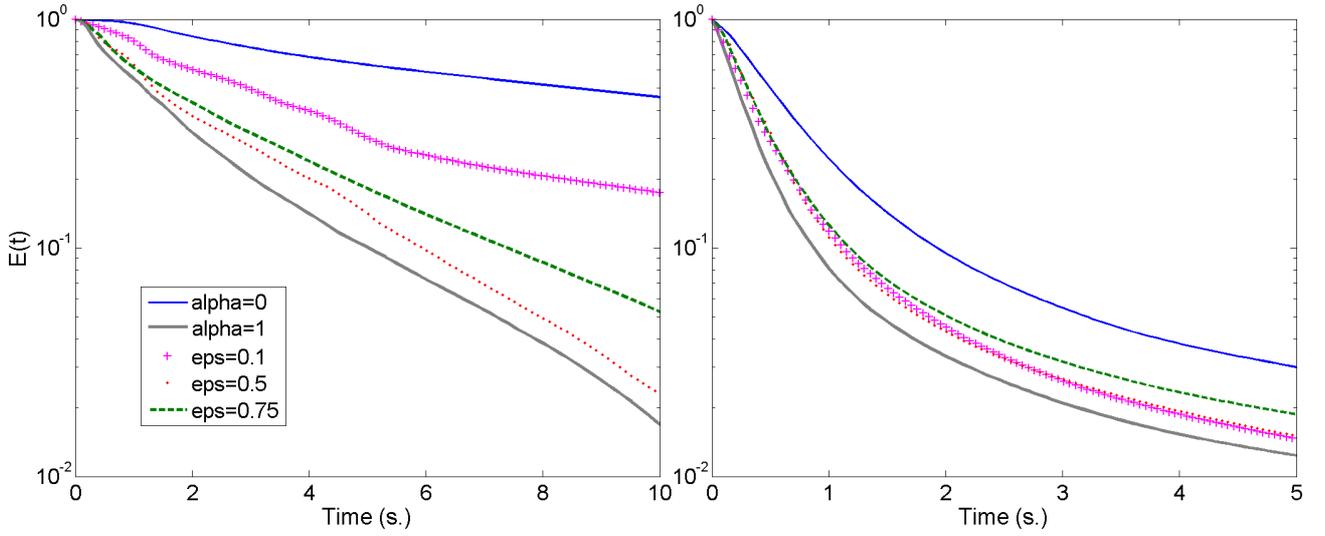}
\caption{Average of functions $E(t)$ for MU using different values of $\epsilon$, with $\alpha=0$ and $\alpha=1$ for reference (see Section~\ref{fnit}): (left) dense matrices, (right) sparse matrices.  Same settings as in Figure~\ref{mualpha}.} 
\label{mueps}
\end{center}
\end{figure*}

\begin{figure*}[ht!]
\begin{center}
\includegraphics[width=\textwidth]{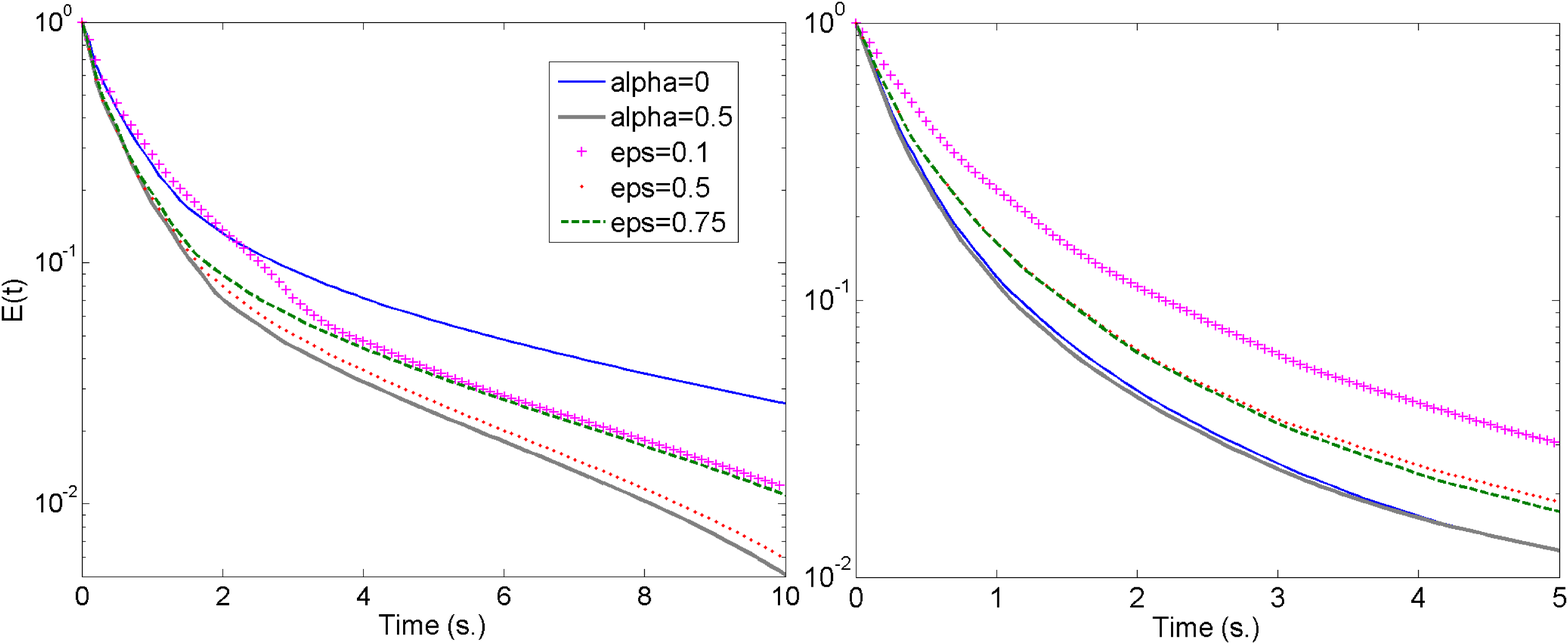}
\caption{Average of functions $E(t)$ for HALS using different values of $\epsilon$, with $\alpha=0$ and $\alpha=0.5$ for reference (see Section~\ref{fnit}): (left) dense matrices, (right) sparse matrices.  Same settings as in Figure~\ref{mualpha}.} 
\label{halseps}
\end{center}
\end{figure*}

Figures~\ref{mueps} shows the results for MU with different values of $\epsilon$ (we also include the original MU and MU with $\alpha = 1$ presented in the previous section to serve as a comparison). 
Figures~\ref{halseps} displays the same experiment for HALS. 

In both cases, we observe that the dynamic stopping criterion is not able to outperform the approach based on a fixed number of inner iterations  ($\alpha = 1$ for MU, $\alpha = 0.5$ for HALS). Moreover, in the experiments for HALS with sparse matrices, it is not even able to compete with the original algorithm. 


\subsection{A Hybrid Stopping Criterion} \label{hybrid}

We have shown in the previous section that using a fixed number of inner iterations works better than a stopping criterion based solely on the difference between two iterates. 
However, in some circumstances, we have observed that inner iterations become ineffective before their maximal count is reached, so that it would in some cases be worth switching earlier to the other factor. 

This occurs in particular when the numbers of rows $m$ and columns $n$ of matrix $M$ have different orders of magnitude. For example, assume without loss of generality that $m \ll n$, so that we have $\rho_W \gg \rho_H$. Hence, on the one hand, matrix $W$ has significantly less entries than $H$ ($mr \ll nr$), and the corresponding NNLS subproblem features a much smaller number of variables; on the other hand, $\rho_W \gg \rho_H$ so that the above choice will lead many more updates of $W$ performed. In other words, many more iterations are performed on the simpler problem, which might be unreasonable.  
For example, for the CBCL face database (cf.\@ Section~\ref{ne}) with $m = 361$, $n = 2429$ and $r = 20$, we have $\rho_H \approx 18$ and $\rho_W \approx 123$, and this large number of inner $W$-updates is typically not necessary to obtain an iterate close to an optimal solution of the corresponding NNLS subproblem. \\

\begin{figure*}[ht!]
\begin{center}
\includegraphics[width=\textwidth]{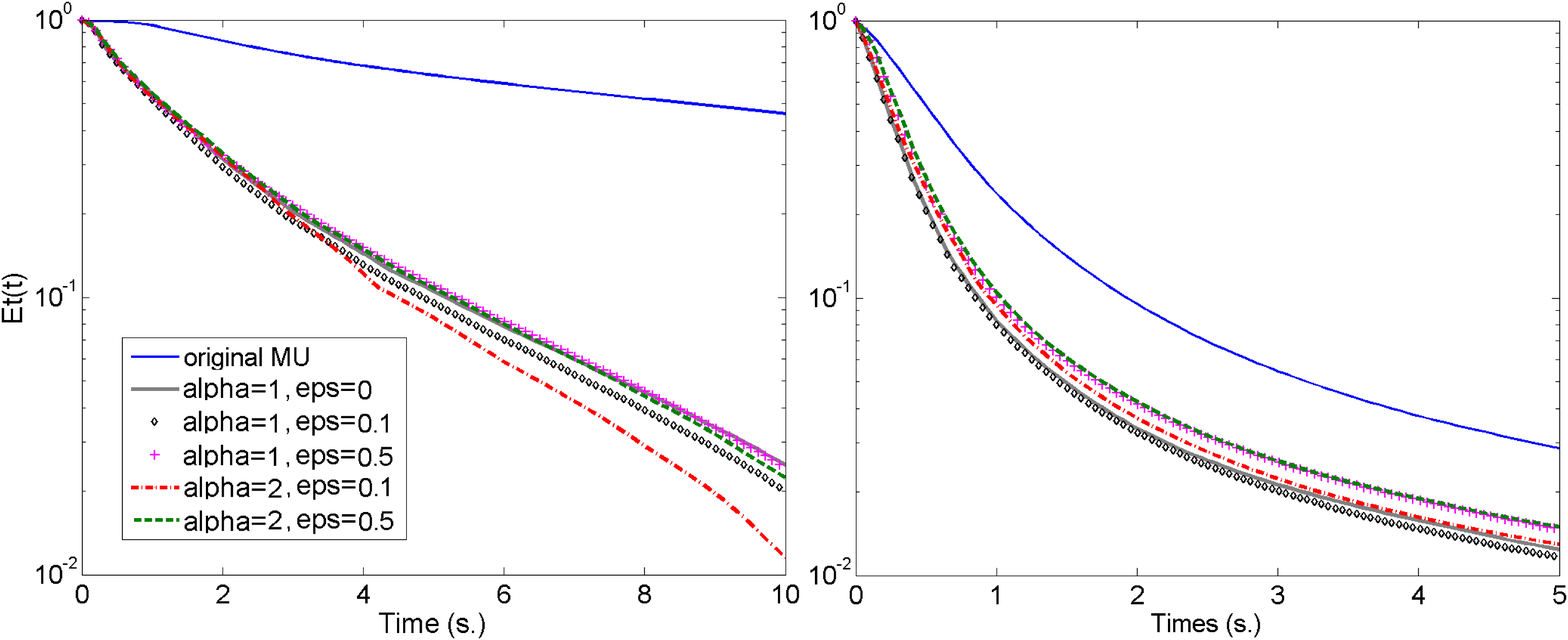}
\caption{Average of functions $E(t)$ for MU using different values of $\alpha$ and $\epsilon$: (left) dense matrices, (right) sparse matrices.  Same settings as in Figure~\ref{mualpha}.} 
\label{mualphaeps}
\end{center}
\end{figure*}

\begin{figure*}[ht!]
\begin{center}
\includegraphics[width=\textwidth]{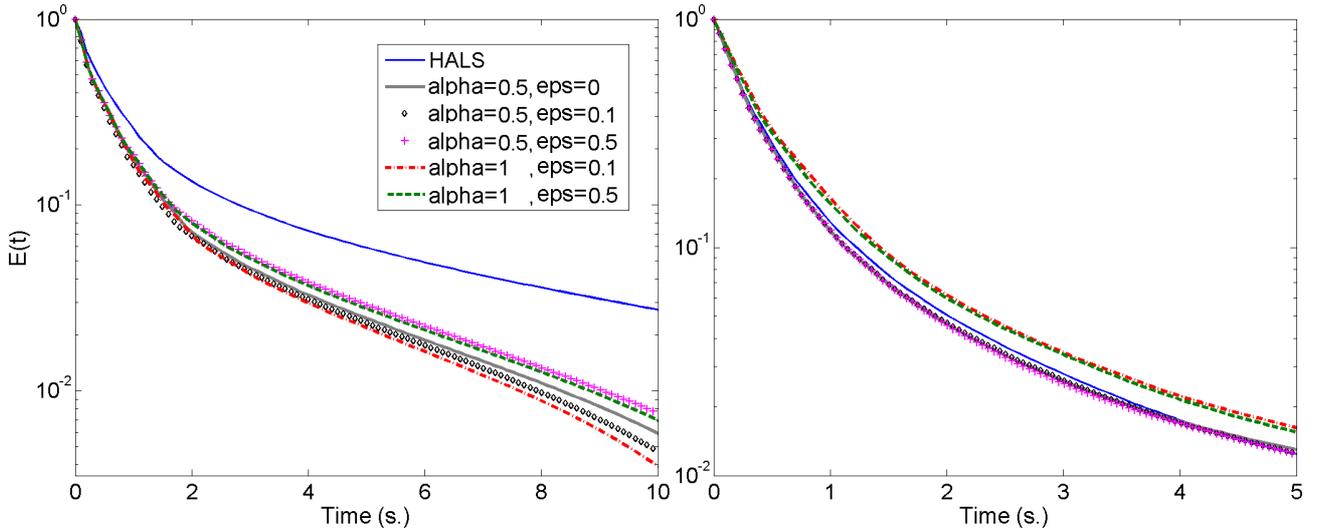}
\caption{Average of functions $E(t)$ for HALS using different values of $\alpha$ and $\epsilon$: (left) dense matrices, (right) sparse matrices.  Same settings as in Figure~\ref{mualpha}.} 
\label{halsalphaeps}
\end{center}
\end{figure*}

Therefore, to avoid unnecessary inner iterations, we propose to combine the fixed number of inner iterations proposed in Section~\ref{fnit} with the supplementary stopping criterion described in Section~\ref{dcit}. 
This safeguard procedure will stop the inner iterations before their maximum number $\lfloor 1+\alpha \rho_W \rfloor$ is reached when they become ineffective (depending on parameter $\epsilon$, see Equation~\eqref{stopcrit}). 
 Algorithm~\ref{anmf} displays the pseudocode for the corresponding accelerated MU, as well as a similar adaptation for HALS. 
\algsetup{indent=2em} 
\begin{algorithm}
\caption{Accelerated MU and HALS} \label{anmf}
\begin{algorithmic}[1]
\REQUIRE Data matrix $M \in \mathbb{R}^{m \times n}_+$ and initial iterates $(W^{(0)},H^{(0)}) \in \mathbb{R}^{m \times r}_+ \times \mathbb{R}^{r \times n}_+$.
\FOR {$k = 0, 1, 2, \dots$} \vspace{0.1cm}
\STATE Compute $A = M {H^{(k)}}^T$ and $B = H^{(k)} {H^{(k)}}^T$;  $W^{(k,0)} = W^{(k)}$;  
 \FOR {$l = 1$ : $\lfloor 1+\alpha \rho_W \rfloor$}
 	\STATE  Compute $W^{(k,l)}$ using either MU or HALS (cf.\@ Algorithms~\ref{MUalgo} and \ref{HALSalgo}); \vspace{0.1cm}    
 	\IF {$||W^{(k,l)}-W^{(k,l-1)}||_F \leq \epsilon ||W^{(k,1)}-W^{(k,0)}||_F$} 
 		\STATE break;
 	\ENDIF
 \ENDFOR
 \STATE  $W^{(k+1)} = W^{(k,l)}$; 
 \STATE Compute $H^{(k+1)}$ from $H^{(k)}$ and $W^{(k+1)}$ using a symmetrically adapted version of steps 2-9; 
 \ENDFOR
\end{algorithmic}
\end{algorithm}
Figures~\ref{mualphaeps} and \ref{halsalphaeps} displays the numerical experiments for MU and HALS respectively. 

In the dense case, this safeguard procedure slightly improves performance. We also note that the best values of parameter $\alpha$ now seem to be higher than in the unsafeguarded case ($\alpha=2$ versus $\alpha=1$ for MU, and $\alpha=1$ versus $\alpha=0.5$ for HALS). Worse performance of those higher values of $\alpha$ in the unsafeguarded scheme can be explained by the fact that additional inner iterations, although sometimes useful, become too costly overall if they are not stopped when becoming ineffective. 

\begin{figure*}[ht!]
\begin{center}
\includegraphics[width=\textwidth/2]{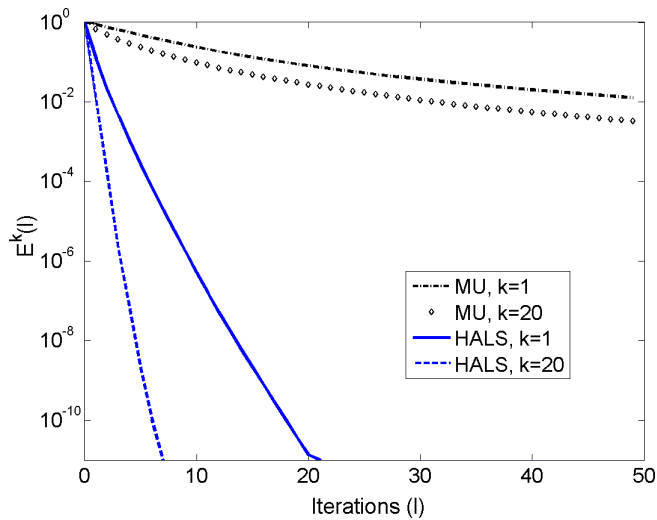}
\caption{Evolution of the relative error $E^k(l)$ of the iterates of inner iterations in MU and HALS, solving the NNLS subproblem $\min_{W \geq 0}||M-WH^{(k)}||_F$ with $r = 40$ for the classic text dataset (cf.\@ Table~\ref{dtm}).} 
\label{denseVSsparse}
\end{center}
\end{figure*}

In the sparse case, the improvement is rather limited (if not absent) and most accelerated variants provide similar performances. In particular, as already observed in Sections~\ref{fnit} and \ref{dcit}, the accelerated variant of HALS does not perform very differently from the original HALS on sparse matrices. We explain this by the fact that HALS applied on sparse matrices is extremely efficient and one inner update already decreases the objective function significantly. To illustrate this, Figure~\ref{denseVSsparse} shows the evolution of the relative error 
\[
E^k(l) = \frac{||M-W^{(k,l)}H^{(k)}||_F-e^k_{\text{min}}}{||M-W^{(k,0)}H^{(k)}||_F-e^k_{\text{min}}}
\]
of the iterate $W^{(k,l)}$ for a sparse matrix $M$, where\footnote{We have used the active-set algorithm of \cite{KP09} to compute the optimal value of the NNLS subproblem.} $e^k_{\text{min}} = \min_{W \geq 0}||M-WH^{(k)}||_F$. Recall that $(W^{(k,0)},H^{(k)})$ denotes the solution obtained after $k$ outer iterations (starting from randomly generated matrices). 
  For $k=1$ (resp.\@ $k=20$), the relative error is reduced by a factor of more than 87\% (resp.\@ 97\%) after only one inner iteration. 
  

\subsection{Application to Lin's Projected Gradient Algorithm} \label{appllin}

The accelerating procedure described in the previous sections can potentially be applied to many other NMF algorithms. To illustrate this, we have modified Lin's projected gradient algorithm (PG) \citep{Lin} by replacing the original dynamic stopping criterion (based on the stationarity conditions) by the hybrid strategy described in Section~\ref{hybrid}. 
It is in fact straightforward to see that our analysis is applicable in this case, since Lin's algorithm also requires the computation of $HH^T$ and $MH^T$ when updating $W$, because the gradient of the objective function in \eqref{NMF} is given by  $\nabla_W ||M-WH||_F^2 = 2 WHH^T - 2 MH^T$. This is also a direct confirmation that our approach can be straightforwardly applied to many more NMF algorithms than those considered in this paper. 

\begin{figure*}[ht!]
\begin{center}
\includegraphics[width=\textwidth]{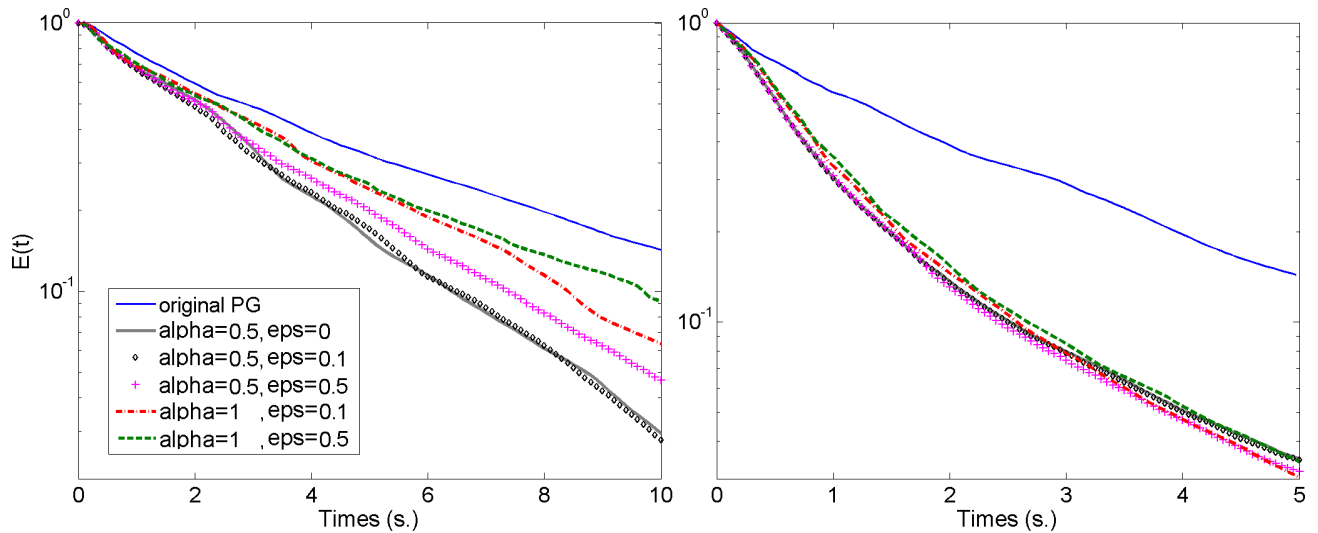}
\caption{Average of functions $E(t)$ for the projected gradient algorithm of \cite{Lin}, and its modification using a fixed number of inner iterations. Same settings as Figure~\ref{mualpha}.} 
\label{linalpha}
\end{center}
\end{figure*}

Figure~\ref{linalpha} displays the corresponding computational results, comparing the original PG algorithm (as available from \cite{Lin}) with its dynamic stopping criterion (based on the norm of the projected gradient) and our variants, based on a (safeguarded) fixed number of inner iterations. It demonstrates that our accelerated schemes perform significantly better, both in the sparse and dense cases (notice that in the sparse case, most accelerated variants perform similarly). The choice $\alpha = 0.5$ gives the best results, and  
 the safeguard procedure does not help much; the reason being that PG converges relatively slowly  (we will see in  Section~\ref{ne} that its accelerated variant converges slower than the accelerated MU). 

\section{Convergence to Stationary Points} \label{conv} 

In this section, we briefly recall convergence properties of both MU and HALS, and show that they are inherited by their accelerated variants. 


\subsection{Multiplicative Updates} \label{muconv}

It was shown by \cite{DM86} and later by \cite{LS1} that a single multiplicative update of $W$ (i.e., replacing $W$ by $W \circ \frac{[MH^T]}{[WHH^T]}$ while $H$ is kept fixed) guarantees that the objective function $||M-WH||_F^2$ does not increase. Since our accelerated variant simply performs several updates of $W$ while $H$ is unchanged (and vice versa), we immediately obtain that the objective function $||M-WH||_F^2$ is non-increasing under the iterations of Algorithm~\ref{anmf}. 


Unfortunately, this property does not guarantee convergence to a stationary point of \eqref{NMF}, and this question on the convergence of the MU seems to be still open, see \cite{Linb}. Furthermore, in practice, rounding errors might set some entries in $W$ or $H$ to zero,  and then multiplicative updates cannot modify their values. Hence, it was observed that despite their monotonicity, MU do not necessarily converge to a stationary point, see \cite{GZ05}. 

However, \cite{Linb} proposed a slight modification of MU in order to obtain the convergence to a stationary point. Roughly speaking, MU is recast as a rescaled gradient descent method and the step length is modified accordingly. 
Another even simpler possibility is proposed by \cite{GG08} who proved the following theorem (see also \cite[\S 4.1]{G11} where the influence of parameter $\delta$ is discussed): 
\begin{theorem}[\cite{GG08}] \label{epsMU} For any constant $\delta > 0$, $M \geq 0$ and any\footnote{$(W, H) \geq \delta$ means that $W$ and $H$ are component-wise larger than $\delta$.} $(W, H) \geq \delta$,  $||M-WH||_F$ is nonincreasing under 
\begin{equation}
W \leftarrow \max\Big(\delta,W \circ \frac{[M H^T]}{[W H H^T]}\Big), \qquad
H \leftarrow \max\Big(\delta,H \circ \frac{[W^T M]}{[W^T W H]} \Big), 
\label{LSupdateEps}
\end{equation}
where the $\max$ is taken component-wise. 
Moreover, every limit point of the corresponding (alternated) algorithm is a stationary point of the following optimization problem
\begin{equation} 
\min_{W\geq \delta, H\geq \delta} ||M-WH||_F^2. \nonumber 
\end{equation}
\end{theorem}
The proof of Theorem~\ref{epsMU} only relies on the fact that the limit points of the updates \eqref{LSupdateEps} are fixed points (there always exists at least one limit point because the objective function is bounded below and non-increasing under updates \eqref{LSupdateEps}). 
Therefore, one can easily check that the proof still holds when a bounded number of inner iterations is performed, i.e., the theorem applies to our accelerated variant (cf.\@ Algorithm~\ref{anmf}). 

\begin{figure*}[ht!]
\begin{center}
\includegraphics[width=\textwidth]{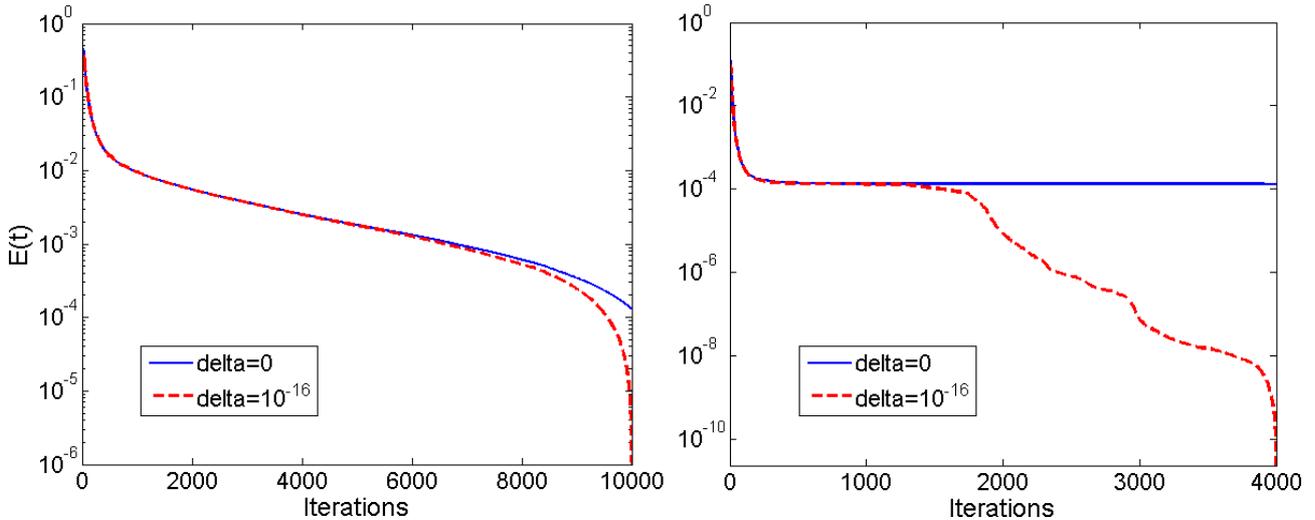} 
\caption{Functions $E(t)$ for $\delta = 0$ and $\delta = 10^{-16}$ on the (dense) ORL face dataset (cf.\@ Table~\ref{dip}) and the (sparse) classic text dataset (cf.\@ Table~\ref{dtm}) with $r=40$.} 
\label{Ceps}
\end{center}
\end{figure*}

It is important to realize that this is not merely a theoretical issue and that this observation can really play a crucial role in practice. To illustrate this, Figure~\ref{Ceps} shows the evolution of the normalized objective function (cf.\@ Equation~\eqref{Evol}) using $\delta = 0$ and $\delta = 10^{-16}$ starting from the same initial matrices $W^{(0)}$ and $H^{(0)}$ randomly generated (each entry uniformly drawn between 0 and 1). 
We observe that, after some number of iterations, the original MU (i.e., with $\delta = 0$) get stuck while the variant with $\delta = 10^{-16}$ is still able to slightly improve $W$ and $H$. Notice that this is especially critical on sparse matrices (see Figure~\ref{Ceps}, right) because many more entries of $W$ and $H$ are expected to be equal to zero at stationarity. For this reason, in this paper, all numerical experiments with MU use the updates from Equation~\eqref{LSupdateEps} with $\delta = 10^{-16}$ (instead of the original version with $\delta = 0$). 


\subsection{Hierarchical Alternating Least Squares} 

HALS is an exact block-coordinate descent method where blocks of variables (columns of $W$ and rows of $H$) are optimized in a cyclic way (first the columns of $W$, then the rows of $H$, etc.).  Clearly, exact block-coordinate descent methods always guarantee the objective function to decrease. However, convergence to a stationary point requires additional assumptions. For example,  \cite{B99, B99b} (Proposition 2.7.1) showed that if the following three conditions hold: 
\begin{itemize}
\item each block of variables belongs to a closed convex set (which is the case here since the blocks of variables belong either to $\mathbb{R}^m_+$ or $\mathbb{R}^n_+$), 

\item the minimum computed at each iteration for a given block of variables is uniquely attained; 

\item the function is monotonically nonincreasing in the interval from one iterate to the next; 
\end{itemize}
then exact block-coordinate descent methods converge to a stationary point. 
The second and the third requirements are satisfied as long as no columns of $W$ and no rows of $H$ become completely equal to zero (subproblems are then strictly convex quadratic programs, whose unique optimal solutions are given by the HALS updates, see Section~\ref{intro}). In practice, if a column of $W$ or a row of $H$ becomes zero\footnote{In practice, this typically only happens if the initial factors are not properly chosen, see \cite{diep}.}, we reinitialize it to a small positive constant (we used $10^{-16}$). 
We refer the reader to \cite{diep} and \cite{GG08} for more details on the convergence issues related to HALS.    

Because our accelerated variant of HALS is just another type of exact block-coordinate descent method (the only difference being that the variables are optimized in a different order: first several times the columns of $W$, then several times the rows of $H$, etc.), it inherits 
all the above properties. In fact, the statement of the above-mentioned theorem in \cite[p.\@6]{B99b} mentions that `the order of the blocks may be arbitrary as long as there is an integer $K$ such that each block-component is iterated at least once in every group of $K$ contiguous iterations', which clearly holds for our accelerated algorithm with a fixed number of inner iterations and its hybrid variant\footnote{Note however that the accelerated algorithms based solely on the dynamic stopping criterion (Section~\ref{dcit}) might not satisfy this requirement, because the number of inner iterations for each outer iteration can in principle grow indefinitely in the course of the algorithm.}.  


\section{Numerical Experiments} \label{ne}

In this section, we compare the following algorithms, choosing for our accelerated MU, HALS and PG  schemes the hybrid stopping criterion and the best compromise for values for parameters $\alpha$ and $\epsilon$ according to tests performed in Section~\ref{accelalgo}. 
\begin{enumerate} 
\item (\textbf{MU}) \, The multiplicative updates algorithm  of \cite{LS2}. 
\item (\textbf{A-MU}) \, The accelerated MU with a safeguarded fixed number of inner iterations using $\alpha = 2$ and $\epsilon = 0.1$ (cf.\@ Algorithm~\ref{anmf}).  
\item (\textbf{HALS}) \, The hierarchical alternating least squares algorithm  of \cite{Cic}.  
\item (\textbf{A-HALS}) \, The accelerated HALS with a safeguarded fixed number of inner iterations using $\alpha = 0.5$ and $\epsilon = 0.1$ (cf.\@ Algorithm~\ref{anmf}). 
\item (\textbf{PG}) \,  The projected gradient method of \cite{Lin}. 
\item (\textbf{A-PG}) \,  The modified projected gradient method of \cite{Lin} using $\alpha = 0.5$ and $\epsilon = 0$ (cf.\@ Section~\ref{appllin}). 
\item (\textbf{ANLS}) \, The alternating nonnegative least squares algorithm\footnote{Code is available at \url{http://www.cc.gatech.edu/~hpark/}.} of \cite{KP09}, which alternatively optimizes $W$ and $H$ exactly using a  block-pivot active set method. Kim and Park showed that their method typically outperforms other tested algorithms (in particular MU and PG)  on synthetic, images and text datasets. 
\end{enumerate} 

All tests were run using \matlab \, 7.1 (R14), on a 3GHz Intel$^{\textrm{\textregistered}}$ Core{\texttrademark}2 dual core processor. We present numerical results on images datasets (dense matrices, Section~\ref{ima}) and on text datasets (sparse matrices, Section~\ref{txt}). Code for all algorithms but ANLS is available at  
\begin{center}
\url{http://sites.google.com/site/nicolasgillis/}. 
\end{center}

\subsection{Dense Matrices - Images Datasets} \label{ima}

\begin{figure*}[ht!]
\begin{center}
\includegraphics[width=\textwidth]{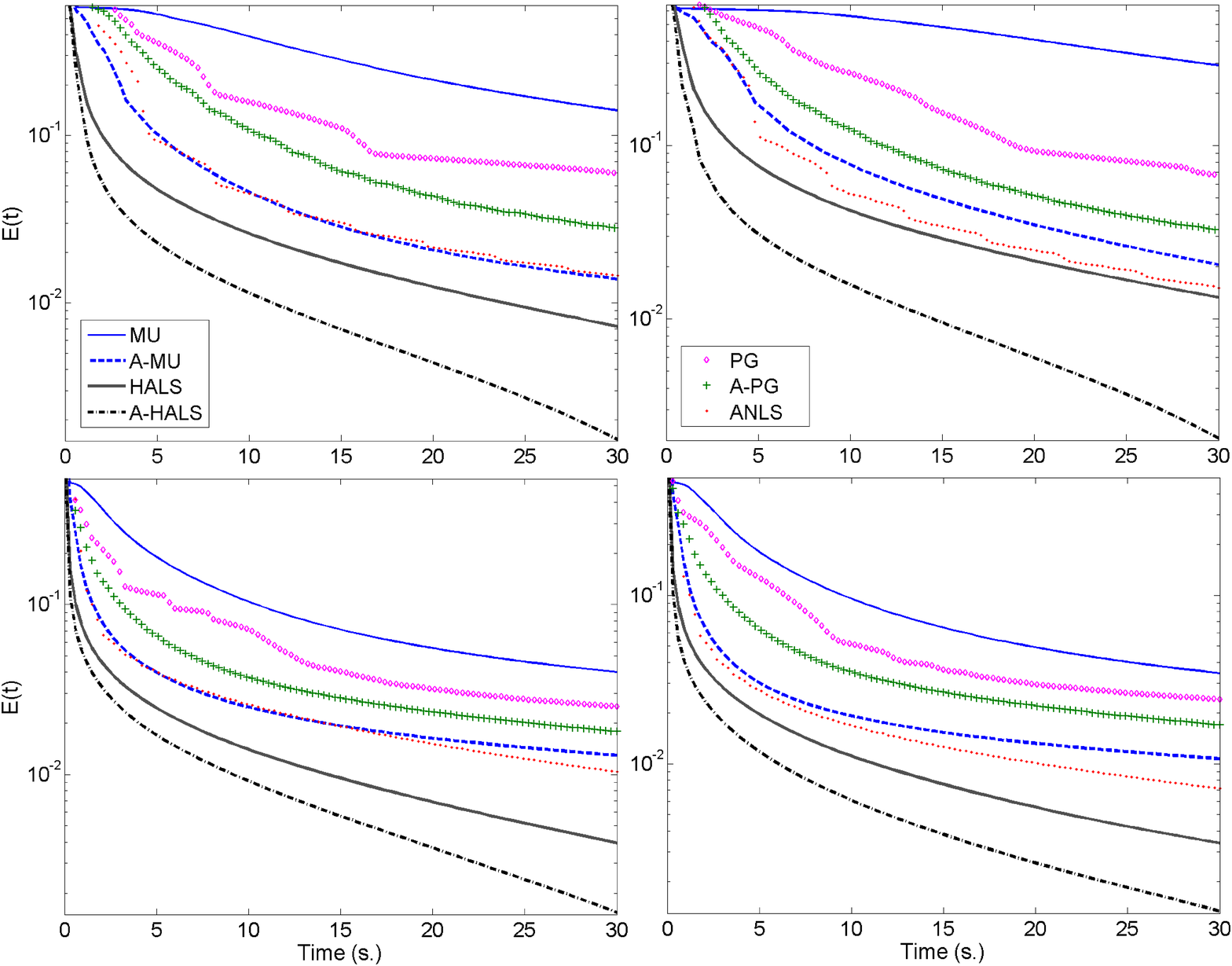} 
\caption{Average of functions $E(t)$ for different image datasets: ORL (top left), Umist (top right), CBCL (bottom left) and Frey (bottom right).} 
\label{densmat}
\end{center}
\end{figure*}

Table~\ref{dip} summarizes the characteristics of the different datasets. 
\begin{center}
\begin{table}[h!]
\begin{center}
\caption{Image datasets. }
\label{dip}  
\begin{tabular}{|c|c|c|c|c|c|c|}
\hline 
Data &             $\#$ pixels &  m  &  n &  r   & $\lfloor  \rho_W \rfloor$ &  $\lfloor  \rho_H \rfloor$ \\ \hline \hline
ORL$^1$    &   $112 \times 92$  & 10304 & 400 & 30, 60  &358, 195 & 13, 7  \\ 
Umist$^2$  &   $112 \times 92$  & 10304 & 575 & 30, 60 & 351, 188 & 19, 10  \\ 
CBCL$^3$ &  $19 \times 19$  & 361  & 2429  & 30, 60 & 12, 7 & 85, 47  \\ 
Frey$^2$  &   $28 \times 20$ & 560 & 1965  & 30, 60 &19, 10 & 67, 36 \\ 
\hline
\end{tabular} \\
\begin{flushleft}

\footnotesize
$\lfloor x \rfloor$ denotes the largest integer smaller than $x$. \\
$^1$ \url{http://www.cl.cam.ac.uk/research/dtg/attarchive/facedatabase.html}\\
$^2$ \url{http://www.cs.toronto.edu/~roweis/data.html}\\
$^3$ \url{http://cbcl.mit.edu/cbcl/software-datasets/FaceData2.html} 
\end{flushleft}
\end{center}
\end{table}
\end{center}
For each dataset, we use two different values for the rank ($r= 30, 60$) and initialize the algorithms with the same 50 random factors $(W^{(0)},H^{(0)})$ (using i.i.d.\@ uniform random variables on $[0,1]$)\footnote{Generating initial matrices $(W^{(0)},H^{(0)})$ randomly typically leads to a very large initial error $e(0)=||M-W^{(0)}H^{(0)}||_F$. 
This implies that $E(t)$ will get very small after one step of any algorithm. To avoid this large initial decrease, we have scaled the initial matrices such that  $\argmin_{\alpha}||M-\alpha W^{(0)}H^{(0)}||_F$=1; this simply amount to multiplying $W$ and $H$ by an appropriate constant, see \cite{GG08}.}. 
In order to assess the performance of the different algorithms, we display individually for each dataset the average over all runs of the function $E(t)$ defined in Equation~\eqref{Evol}, see Figure~\ref{densmat}. 

First, these results confirm what was already observed by previous works: PG performs better than MU \citep{Lin}, ANLS performs better than MU and PG \citep{KP09}, and HALS performs the best \citep{diep}. 
Second, they confirm that the accelerated algorithms indeed are more efficient: A-MU (resp.\@ A-PG) clearly outperforms MU (resp.\@ PG) in all cases, while A-HALS is, by far, the most efficient algorithm for the tested databases. 
It is interesting to notice that A-MU performs better than A-PG, and only slightly worse than ANLS, often decreasing the error as fast during the first iterations. 

\subsection{Sparse Matrices - Text Datasets} \label{txt}

\begin{figure*}[ht!]
\begin{center}
\includegraphics[width=\textwidth]{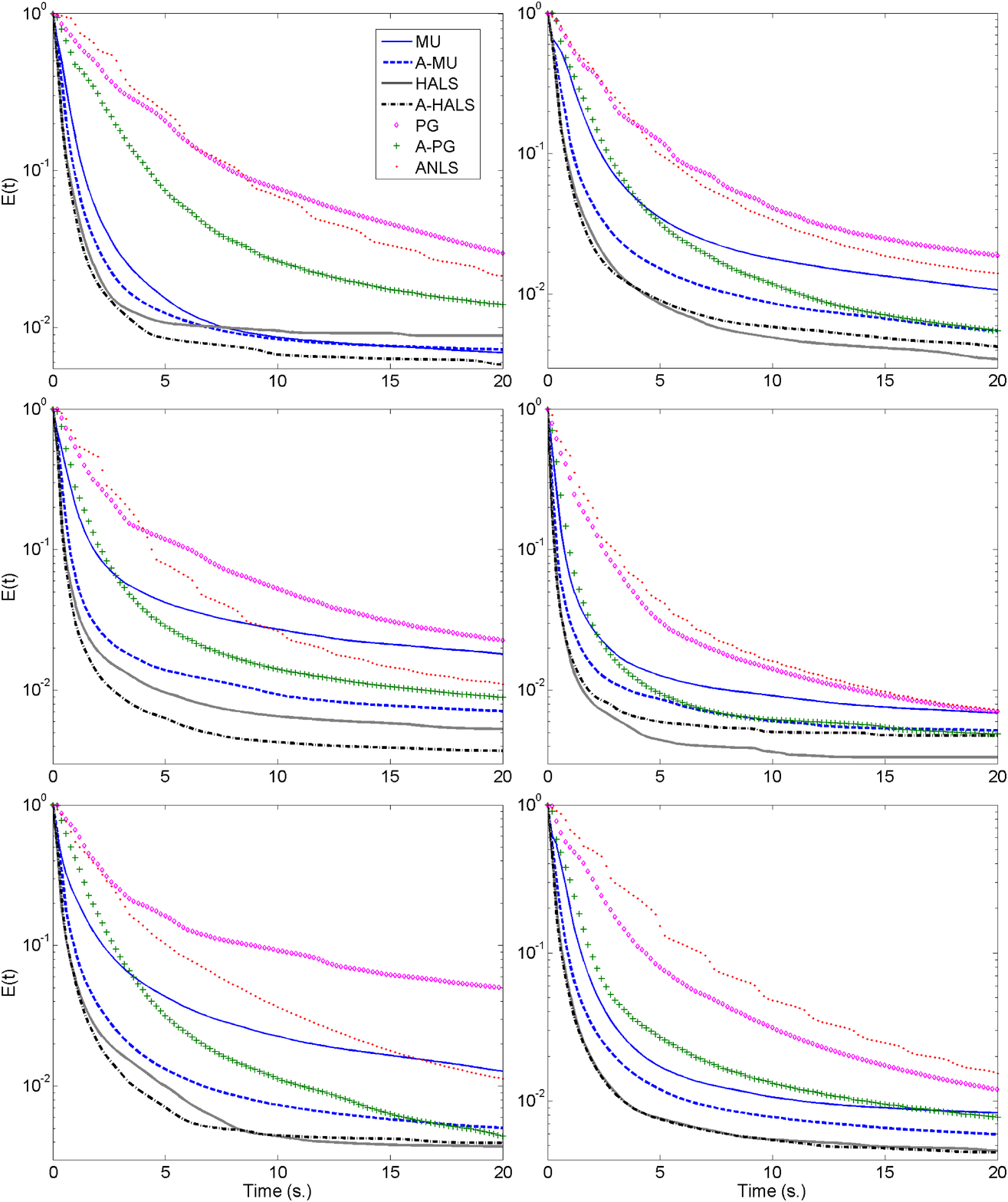}
\caption{Average of functions $E(t)$ for text datasets: classic (top left), sports (top right), reviews (middle left), hitech (middle right), ohscal (bottom left) and la1 (bottom right).} 
\label{sparsemat}
\end{center}
\end{figure*}

Table~\ref{dtm} summarizes the characteristics of the different datasets. 
\begin{center}
\begin{table}[h!] 
\begin{center}
\caption{Text mining datasets \citep{ZG05} (sparsity is given in $\%$: $100*\#\text{zeros}/(mn)$).} 
\label{dtm}
\begin{tabular}{|c|c|c|c|c|c|c|c|} 
\hline
Data &   $m$  &  $n$ &  r &  $\#\text{nonzero}$ & sparsity  & $\lfloor  \rho_W \rfloor$ &  $\lfloor  \rho_H \rfloor$ \\ \hline \hline
classic &  7094   & 41681 & 10, 20 & 223839  & 99.92 & 12, 9 & 2, 1\\ 
sports &    8580 & 14870 &  10, 20 & 1091723 & 99.14 &  18, 11 & 10, 6\\ 
reviews &   4069  & 18483 & 10, 20 & 758635  & 98.99 & 35, 22 & 8, 4\\ 
hitech &    2301 & 10080 &  10, 20 & 331373 & 98.57&  25, 16 & 5, 4\\ 
ohscal &   11162  & 11465 & 10, 20 & 674365 & 99.47 & 7, 4 & 7, 4 \\ 
la1 &    3204 & 31472 &  10, 20 & 484024  & 99.52 & 31, 21 & 3, 2 \\ 
\hline
\end{tabular}
\end{center}
\end{table}
\end{center} 
The factorization rank $r$ was set to 10 and 20.  For the comparison, we used the same settings as for the dense matrices. Figure~\ref{sparsemat} displays for each dataset the evolution of the average of functions $E(t)$ over all runs.   
Again the accelerated algorithms are much more efficient. In particular, A-MU and A-PG converge initially much faster than ANLS, and also obtain better final solutions\footnote{We also observe that ANLS no longer outperforms the original MU and PG algorithms, and only sometimes generates better solutions.}. 
A-MU, HALS and A-HALS have the fastest initial convergence rates, and HALS and A-HALS generate the best solutions in all cases. 
Notice that A-HALS does not always obtain better final solutions than HALS (still this happens on half of the datasets), because HALS already performs remarkably well (see discussion at the end of Section~\ref{hybrid}). However, the initial convergence of A-HALS is in all cases at least as fast as that of HALS.

\section{Conclusion} 

In this paper, we considered the multiplicative updates of \cite{LS2} and the hierarchical alternating least squares algorithm of \cite{Cic}. We introduced accelerated variants of these two schemes, based on a careful analysis of the computational cost spent at each iteration, and preserve the convergence properties of the original algorithms. The idea behind our approach is based on taking better advantage of the most expensive part of the algorithms, by repeating a (safeguarded) fixed number of times the cheaper part of the iterations. This technique can in principle be applied to most NMF algorithms; in particular, we showed how it can substantially improve the projected gradient method of \cite{Lin}. 
We then experimentally showed that these accelerated variants, despite the relative simplicity of the modification, 
significantly outperform the original ones, especially on dense matrices, and compete favorably with a state-of-the-art algorithm, namely the ANLS method of \cite{KP09}. A direction for future research would be to choose the number of inner iterations in a more sophisticated way, with the hope of further improving the efficiency of A-MU, A-PG and A-HALS. 
Finally, we observed that HALS and its accelerated version are the most efficient variants for solving NMF problems, sometimes by far. 


\small

\subsection*{Acknowledgments}

We would like to thank the anonymous reviewers for their insightful comments, which helped to improve the paper.


\end{document}